\documentstyle[twoside,amssymb,12pt]{article}
\newtheorem{prop}{Proposition}[section]

\newtheorem{theo}[prop]{Theorem}

\newtheorem{coro}[prop]{Corollary}

\title{\sc On the case of Kovalevskaya and new examples of  integrable conservative
systems on $S^2$ } 

\author{{\sc K.P. Hadeler}  \\ 
\small  { Mathematische Fakult\"at, Universit\"at T\"ubingen}   \\
\small  { Auf der Morgenstelle 10,  72076  T\"ubingen, Germany}  \\
\small  { e-mail: \ \ hadeler@uni-tuebingen.de }
\and 
{\sc E.N. Selivanova}\thanks{Supported by DAAD.}\\
\small  { Department of Geometry, Nizhny Novgorod State 
Pedagogical University}\\
\small  { 603000 Russia, Nizhny Novgorod, ul. Ulyanova 1}\\
\small  {e-mail: \ \ lena@moebius.mathematik.uni-tuebingen.de} 
}

\begin{document}
\date{}

\maketitle

\noindent
{\bf Abstract.} {\small There  is a well-known example of an integrable 
conservative system on $S^2$,  the case of Kovalevskaya 
in the dynamics of a rigid body, possessing an integral
 of fourth degree in momenta. The aim of this paper is 
to construct  new families of examples of conservative systems on $S^2$
 possessing an 
integral  of  fourth degree in momenta.}

\thispagestyle{empty}

\section{\bf Introduction}

\noindent
Let $M$ be an $n$-dimensional  Riemannian manifold, and $U: M\to \bf R$ be a 
smooth function on $M$. For the Lagrangian $L: 
TM \to \bf R$
we choose $$L(\eta)=\frac{|\eta|^2}{2} + (U\circ\pi )(\eta)$$
where $\pi :TM\to M$ is the canonical projection. In local coordinates $q_1, \ ...\ ,  q_n$, $\dot q_1, \
...\ ,  \dot q_n$  on $TM$ we have
$$L(\eta)=\frac{1}{2}\sum {g_{ij}\dot q^i\dot q^j} + U(q).$$
Identifying $TM$ and $T^*M$ by means of the Riemannian metric, 
we get a Hamiltonian system
with the Hamiltonian $H: T^*M \to \bf R$ which 
 in local coordinates $q_1, \ ...\ , q_n$, $p_1, \ ... \ , p_n$ on $T^*M$ has the form
$$H=\frac{1}{2}\sum{g^{ij}p_ip_j} + U(q)=K+U.$$
We will call these Hamiltonian systems  {\it conservative systems} 
on $M$.

A smooth function $F : T^*M \to \bf R$ which is an integral of the 
Hamiltonian system with the Hamiltonian 
$H$ and which is independent of $H$ we will call an {\it integral} of this system 
of {\it degree} $m$ in momenta if in local coordinates $F$ has the form
$$F=\sum_{k_1+ \ ...\  + k_n\le m}{a_{k_1 \ ...\    k_n}(q)p_1^{k_1}  \ ...\   p_n^{k_n}}.$$ 

We will say that two Hamiltonians $H_1=K_1+U_1$ and $H_2=K_2+U_2$ 
are {\it equivalent} if there exists a diffeomorphism  $\phi$ of $M$ and a 
diffeomorphism
$\Phi$ of  $T^*M$
such that the diagram
$$\Phi: T^*M \to T^*M$$
$$\pi' \downarrow \ \ \ \ \ \ \ \ \ \ \downarrow \pi' $$
$$\phi: M\to M$$
is commutative, where   $\Phi$ is linear for   $p\in M$ fixed, and if
there are  some nonzero constants $\kappa$, $\tilde \kappa$ such that
$\Phi^*(K_1)=\kappa K_2$, $\phi^*(U_1)=\tilde \kappa U_2$.

Clearly,  if the Hamiltonians $H_1$, $H_2$ are equivalent
and one of the corresponding systems possesses an integral of degree $m$ in momenta,
then the other system has the same property.

\medbreak
There  is a well-known example of an integrable 
conservative system on $S^2$,  the case of Kovalevskaya \cite{Kov} 
in the dynamics of a rigid body, possessing an integral
 of fourth degree in momenta. 
The total energy (the Hamiltonian) of this 
system has the following form
\begin{equation}
H=\frac{du_1^2 + du_2^2 + 2du_3^2}{2u_1^2 + 2u_2^2 + u_3^2} - u_1,
\label{kov1}
\end{equation}
where $S^2$ is given by $ u_1^2 + u_2^2 + u_3^2=1$, (see 
 \cite{BKF}).

Goryachev  proposed in \cite{Gor} a  
family of examples of conservative systems on $S^2$
possessing an integral of fourth
degree in momenta. The Hamiltonians of these systems have the following form
\begin{equation}
H=\frac{du_1^2 + du_2^2 + 2du_3^2}{2u_1^2 + 2u_2^2 + u_3^2} - 2B_1u_1u_2 - B_2(u_1^2-u_2^2) -u_1
\label{gor}
\end{equation}
where $S^2$ is given by $ u_1^2 + u_2^2 + u_3^2=1$, $B_1$, $B_2$ are arbitrary constants.
Clearly, this family reduces to 
the case of Kovalevskaya when $B_1=B_2=0$.

\medbreak
The aim of this paper is to construct {\it new families of examples} of conservative systems on $S^2$
 possessing an 
integral $F$  of {\it fourth degree} in momenta, i.e. examples which are not 
equivalent to the cases of Kovalevskaya or Goryachev.

 We will show that these families include the integrable cases
from \cite{Se} and the case of Kovalevskaya. In particular, we will obtain
an explicit expression for the case of Kovalevskaya from a quite different point
of view.

\section{New examples }

In \cite{Se} the following local criterion for integrability 
of the geodesic flow of a Riemannian metric $ds^2$ by a polynomial
of fourth degree in momenta has been obtained.

\begin{theo} Let $ds^2=\lambda(z, \bar z)dzd\bar z$ be a metric
such that there exists a function $f:{\bf R}^2\mapsto \bf R$, satisfying the 
following conditions 
\begin{equation}
\lambda(z,\bar z)=\frac{\partial^2f}{\partial  z\partial \bar z}, \ \ 
\mbox{Im}\ \left(\frac{\partial^4f}{\partial  z^4}
\frac{\partial^2f}{\partial  z\partial \bar z} + 
3\frac{\partial^3f}{\partial  z^3}\frac{\partial^3f}
{\partial  z^2\partial \bar z} + 2\frac{\partial^2f}{\partial  z^2}
\frac{\partial^4f}{\partial  z^3\partial \bar z}\right)=0.
\label{eqpde}
\end{equation}
Then the geodesic flow of $ds^2$ possesses an integral of fourth degree in momenta.
\medbreak
If the geodesic flow of a metric $ds^2$ possesses an integral 
which is a polynomial 
 of  fourth degree in momenta and it does not depend on 
the Hamiltonian and an  integral of smaller degree then there exist
  conformal coordinates $x, y$ and a function   $f:{\bf R}^2\mapsto \bf R$ such that 
$ds^2=\lambda(z, \bar z)dzd\bar z$ where $z=\varphi +i y$ and (\ref{eqpde}) holds.
\label{loc-th}
\end{theo} 

Equation (\ref{eqpde}) in some other form has been obtained also in 
\cite{Hall}.

As in \cite{Se} we consider now the solutions of (\ref{eqpde}) 
of the following form:
$$
f(\varphi,y)=u(y)\cos \varphi + \xi(y) + d(\varphi^2 - y^2)
$$
where $u(y)$, $\xi(y)$ are some smooth functions and $d$ is a constant.
In \cite{Se} it has been shown that the geodesic flows of the metrics
$ds^2=\lambda(d\varphi^2+dy^2)$ where
$$
\lambda=\frac{1}{4}\left((u''(y)- u(y))\cos \varphi + \xi''(y) \right),
$$

$$\xi''=\frac{d_1u(y) + c}{(u'(y))^2}, \ d_1, c - const, \quad \mbox{for} \quad  d=0, $$
 or 
$$\xi''=2d\frac{u'^2(y) - u^2(y) +
 d_1(2d)^{-1}u(y) + p}{u'(y))^2}, \  d_1,  p - const, \quad \mbox{for} \quad d\ne 0, $$
and
$u(y)$ satisfies
\begin{equation}
2u''^2 - 3u^2 + u'u''' = \frac{a}{2},\ a - const,
\label{0}
\end{equation}
possess an integral of fourth degree in momenta.

Using the   well-known Maupertuis's principle and taking into account that 
$c$ and $d$ are  arbitrary constants we find  that
the Hamiltonian systems with the Hamiltonians
\begin{equation}
H=\frac{d\varphi^2 + dy^2}{u'^2(y)} -
 (u''(y)-u(y))u'^2(y)\cos \varphi
\label{H}
\end{equation}
and 
\begin{equation}
H_p=\frac{u'^2(y)-u^2(y)+p}{u'^2(y)}\left(d\varphi^2 + dy^2\right) -
 (u''(y)-u(y))\frac{u'^2(y)}{u'^2(y)-u^2(y)+p}\cos \varphi
\label{H_p}
\end{equation}
where   $u$ is a solution of (\ref{0}),
possess also an integral of fourth degree in momenta.

Now we exploit the properties of the differential equation (\ref{0}).
\begin{prop}
The differential equation (\ref{0}) is equivalent to the two-parameter
family of first-order differential equations
\begin{equation}
u'^4=b+b_1u + au^2 + u^4
\label{1od}
\end{equation}
where $b$, $b_1$ are arbitrary constants.
\label{sol}
\end{prop}

\noindent
{\em Proof.} Multiply equation (\ref{0}) by $u'$
$$ u'''u'^2 + 2u''u'u'' - 3u^2u'= \frac{a}{2}u'$$
and integrate
$$u''u'^2 - u^3 =  \frac{a}{2}u + \frac{b_1}{4} $$
with some constant $b_1$.
Multiply again with $u'$
$$u''u'^3 - u^3u'= \frac{a}{2}uu' + \frac{b_1}{4}u'$$
and integrate
$$\frac{1}{4}u'^4 - \frac{1}{4}u^4=\frac{a}{4}u^2 + \frac{b_1}{4}u + \frac{b}{4}$$
with some constant $b$. This is (\ref{1od}).

\hfill$\Box$

Thus the family (\ref{H}) can be parametrized by $a$, $b$, $b_1$, $u(0)$  and (\ref{H_p}) by $a$, $b$, $b_1$, $u(0)$, 
$p$.
\medbreak
We will consider the case $b_1=0$ in (\ref{1od}).
Denote $ {A}(u,a,b)=b+au^2+u^4$ where $a$ and $b$ are 
constants. 

Thus, the construction of our examples is based on the properties
of the following differential equation
\begin{equation}
u'^4=b+au^2+u^4.
\label{2od}
\end{equation}
Due to Proposition \ref{sol} we may obtain  explicit expressions for 
 the Hamiltonians (\ref{H})
and (\ref{H_p}) in the coordinates $\varphi$ and $u$:
\begin{equation}
H= {A}^{-\frac{1}{2}}(u,a,b)
\left(d\varphi^2 +  {A}^{-\frac{1}{2}}du^2\right) - 
\frac{u}{2}((a+2u^2)-2 {A}^{\frac{1}{2}})\cos \varphi
\label{H1}
\end{equation}
and
\begin{equation}
H_{p}=\frac{ {A}^{\frac{1}{2}}-u^2+p}{ {A}^{\frac{1}{2}}}
\left(d\varphi^2 +  {A}^{-\frac{1}{2}}du^2\right) - 
\frac{u((a+2u^2) -2{A}^{\frac{1}{2}})}{2( {A}^{\frac{1}{2}}-u^2+p)}
\cos \varphi
\label{H_p1}
\end{equation}

Thus, in $(u, \varphi)$ coordinates, $H$ depends only  the parameters $a$, $b$, and 
$H_p$ on $a$, $b$, $p$.

By a substitution $u\to \alpha_1 u$, $\alpha_1 - const$ in (\ref{2od}) we obtain
$$u'^4=\frac{b}{\alpha_1^4} + \frac{a}{\alpha_1^2}u^2 + u^4.$$
So we can normalize $b$ to $1$, $0$ or $-1$ and keep $a$.
\medbreak
In this section we will consider the case $b=1$. 
In the next section
we will consider the case $b=0$, $a=1$ and show that it is in fact the 
case of Kovalevskaya.
\medbreak
The case $b=1$, $a=0$ has been considered in \cite{Se}. It has been shown that
in this case the systems given by (\ref{H}), (\ref{H_p})
(and therefore (\ref{H1}), (\ref{H_p1})) define smooth conservative systems
on $S^2$ possessing an integral of fourth degree in momenta.
\bigbreak
We will use the following proposition.

\begin{prop}
If $b=1$ in (\ref{2od}) and $a>-2$, then 
there is a solution $u_a(y): {\bf R}\to {\bf R}$ of (\ref{2od}) such that
the following holds
$$u_a'(y)=(\exp y)\nu_a(\exp(-2y))=\exp(- y)\nu_a(\exp(2y)),$$
$$u_a'^2(y)(u_a''(y) - u_a(y))=\exp (-y)\mu_a(\exp(-2y))=-
(\exp y)\mu_a(\exp(2y))$$
where the functions $\mu_a$, $\nu_a$ are of class $C^{\infty}$
and  $\nu_a>0$ everywhere. 
\label{cen}
\end{prop}

\noindent
{\em Proof.} 
By a simple computation we get that
if $-2<a<0$, then $A(x, a, 1)\ge 1- \frac{a^2}{4}>0$ for all $x\in \bf R$
and if $a\ge 0$, then $A(x, a, 1)\ge 1$ for all $x\in \bf R$.

Therefore, if $b=1$ and $a>-2$, then all solutions $u$ of (\ref{2od})  exist globally.
All increasing (descreasing) solutions are translates  of each other. Increasing
  solutions run from $-\infty$ to $ +\infty$, similarly for descreasing solutions.

W.l.o.g. we may only consider increasing solutions of
\begin{equation}
u'=(1+au^2 + u^4)^{\frac{1}{4}}, \quad u'>0.
\label{new}
\end{equation}

Then we have for $u\ge 0$
$$u'\le u + c_0$$
with a constant $c_0$ such that $c_0\ge 1$ and $6c_0^2\ge a$ and therefore
$$u(y)\le (u(0) + c_0)(\exp y) - c_0$$
and
\begin{equation}
\frac{u(y)}{\exp y}\le u(0) + c_0.
\label{c'}
\end{equation}
Now consider any such solution. Put $s=\exp(-2y)$ and 
$$g(s)=\sqrt{s}u\left(-\frac{1}{2}\log s\right).$$
Then $u(y)=(\exp y)g(\exp(-2y))$. The function $g$ is of class $C^{\infty}$ on $(0, \infty)$.
We normalize $u(0)=g(1)=0$.

Equation (\ref{new}) can be rewritten as a differential 
equation for $g(s)$, $s=\exp(-2y)$:
$$g'=-\frac{s+ag^2}{2((s^2+ag^2 + g^4)^{\frac{1}{4}}+ g)((s^2+ag^2 + g^4)^{\frac{1}{2}}+ g^2)}=\Theta_a(s, g).$$
Thus,  $g$ is decreasing, $g(0)$ is finite in view of (\ref{c'}), and
$g(0)>g(1)=0$.

Therefore, there is a solution $u_a$ of (\ref{2od}) which can be given as
$$u_a(y)=(\exp y)g_a(\exp (-2y))=-(\exp (-y))g_a(\exp 2y)$$
where $g_a$ is of class $C^{\infty}$ on $[0, +\infty)$ if $a\ne 0$.

The case $a=0$ has been considered in \cite{Se} but we can give here another proof.
We can consider the function $\beta(s^2)=g_0(s)$ and prove, with the same arguments 
as above,  that $\beta$ is smooth in zero, and therefore, $g_0$ is 
smooth in zero, too.

Now the corresponding expressions for the functions $\nu_a$ and $\mu_a$
can be obtained in terms of $g_a$.

\hfill$\Box$

Further we will use 
 some properties of the geodesic flows of metrics
\begin{equation}
ds^2=\lambda(r^2)(r^2d\varphi^2 + dr^2)
\label{polya}
\end{equation}
on $S^2$ which have been proved in \cite{Se}.
These properties follow also
from the results of Kolokol'tsov, 
published in his Ph.D. Dissertation,
(Moscow State University, 1984).

\begin{prop}
The geodesic flow of a Riemannian
metric (\ref{polya}) on $S^2$
does not possess a nontrivial 
integral quadratic in momenta (which does not
depend on $H$ and  linear integrals).
\label{kva}
\end{prop}

It is known that a metric of constant positive curvature
has the following
form in polar coordinates
\begin{equation}
ds^2=\frac{C_1}{(1+Dr^2)^2}\left(r^2d\varphi^2 +dr^2\right), 
\label{ccur}
\end{equation}
where $C_1, D - const$.

\begin{prop} 
The geodesic flow of a metric (\ref{polya}) on $S^2$
possesses two independent linear  integrals if and only if it is
 has the form (\ref{ccur}), i.e. if it is 
a metric of constant positive curvature. 
\label{kva2}
\end{prop}

\begin{coro}
Liouville  coordinates $\varphi, y=\log r$, related
to polar coordinates $\varphi, r$
of a metric (\ref{polya}) 
on $S^2$ are unique up to shifts and the transform $y\to -y$.
\label{typ}
\end{coro}

Now we will prove the main theorems.

\begin{theo} Assume 
$b=1$. Then for any  $a>-2$ the Hamiltonian (\ref{H1}) where 
$\varphi\in [0, 2\pi)$ and $u\in (-\infty, +\infty)$
defines a  conservative system
on $S^2$ possessing an integral of fourth degree in momenta.

If $a\ne 2$, then this integral is nontrivial, i.e. 
there is no quadratic or linear integral. 
If $a=2$ it is a Hamiltonian of a metric of constant positive curvature.

The  Hamiltonians (\ref{H1}) for different values of the parameter $a$
are not equivalent.
\label{th1}
\end{theo}

\noindent
{\em Proof.} The integrability of these systems 
follows immediately from Theorem 
\ref{loc-th}, see above. 

So, we have to prove only that the Hamiltonian (\ref{H1}),
for $b=1$ and $a>-2$,  is a sum of a smooth Riemannian metric on $S^2$ and 
a smooth function on $S^2$.

Using Proposition \ref{cen} we may rewrite the corresponding Hamiltonian
in polar coordinates $\varphi, r=\exp y$. By computation we obtain 
$$H=\frac{1}{\nu_a^2(r^2)}(r^2d\varphi^2 + dr^2) + \mu_a(r^2)r\cos \varphi$$
$$=\frac{1}{\nu_a^2(\tilde r^2)}
(\tilde r^2d\tilde \varphi^2 + d\tilde r^2) - 
\mu_a(\tilde r^2)\tilde r\cos \tilde \varphi $$
where $\tilde r=\frac{1}{r}$, $\tilde \varphi=-\varphi$.
Since $\nu_a$, $\mu_a$ are of class $C^{\infty}$ and $\nu_a\ne 0$, 
see Proposition \ref{cen},  this system is a  
conservative system on $S^2$.

\medbreak
Now we prove that the systems with  Hamiltonians (\ref{H1})
where $b=1$ and $a>-2$ do not have linear or  nontrivial quadratic  integrals.

Write $H=\hat H+V(\varphi, y)$ where $\hat H$ is the Hamiltonian of the geodesic 
flow of the metric $$ds_1^2=\frac{d\varphi^2 + dy^2}{u'^2(y)}$$
and 
\begin{equation}
V(\varphi, y)=-u'^2(y)(u''(y)-u(y))\cos \varphi 
\label{V}
\end{equation}
where $u(y)$ satisfies (\ref{2od}) with the parameters $b=1$, $a>-2$.
Note that $ds_1^2$ has the form (\ref{polya}) in $\varphi, r=\exp y$ but it has 
the form (\ref{ccur}) if and only if $a=2$ (if $b=1$). So, we can apply Propositions 
\ref{kva} and \ref{kva2}. We conclude that if $a\ne 2$, then 
 an integral quadratic in momenta of the geodesic flow of $ds_1^2$ depends on the linear
integral $p_{\varphi}$ and the Hamiltonian $\hat H$.

Let us assume that  a system from our theorem has an integral
 which is independent
of  the energy $H$ and which is quadratic in momenta (clearly, this assumption includes
 the case of linear integrals). 

So, there is an  integral $\tilde F$ of (\ref{H1}) which 
is quadratic in momenta.
Thus, $\tilde F= D(p_{\varphi}, p_{y}, \varphi, y) + B(\varphi, y)$ 
where $ D(p_{\varphi}, p_{y}, \varphi, y)$
is a polynomial of  second degree in momenta $p_{\varphi}$, $p_{y}$.

We write $\{\tilde F, H\}=\{
D(p_{\varphi}, p_{y}, \varphi, y) + B(\varphi, y), 
 \hat H +V\}\equiv 0$
and, therefore, $\{D(p_{\varphi}, p_{y}, \varphi, y),  \hat H\}\equiv 0$. 

Thus, the geodesic flow of $ds_1^2$ has an integral quadratic in momenta.
As mentioned above, this integral depends on $p_{\varphi}$ and $\hat H$. Since
$\tilde F$ does not depend on $H$,  we may put w.l.o.g. that
$$D(p_{\varphi}, p_{y}, \varphi, y)=p_{\varphi}^2.$$

We write now
$$\{D,V\}+\{B, \hat H\}=\{p_{\varphi}^2,V\}+\{B,  \hat H\}\equiv 0.$$
By computation we obtain
$$\frac{\partial B}{\partial y}\equiv 0,$$ and
$$u'^{2}(y)\frac{\partial B}{\partial \varphi}=\frac{\partial V}
{\partial \varphi}.$$
So, we get
$$V=(B(\varphi)+\alpha(y))u'^{2}(y)$$
for a smooth function $\alpha(y)$. 
Comparing now this expression with (\ref{V}), we  
 get 
$u''(y)-u(y)\equiv const$, that is  not true if $a\ne 2$.
So, there is no nontrivial quadratic integral of the system given by (\ref{H1})
where $b=1$ and $-2<a<2$ or $a>2$.

\medbreak
In order to prove that the   Hamiltonians (\ref{H1}) for different values of the parameter $a$
are not equivalent
we will apply Corollary \ref{typ}. 

Suppose that the  Hamiltonians (\ref{H1}) for $a=a_1$ and $a=a_2$, ($a_1\ne a_2$) 
are equivalent.
Therefore, from Corollary \ref{typ} 
the following holds identically
$$\frac{du}{d\tilde u}=\pm\left(\frac{1+a_1u^2+u^4}{1+a_2\tilde u^2+\tilde u^4}
 \right)^{\frac{1}{4}}$$ and
$$\frac{1}{\sqrt{1+a_1u^2+u^4}}=\kappa\frac{1}{\sqrt{1+a_2\tilde u^2+\tilde u^4}}$$
for some constant $\kappa$. Thus, we obtain $u=\tilde u$ and, therefore,
 $a_1=a_2$.

\medbreak
With the same arguments as in \cite{Se}, we may prove that no Hamiltonian  from our theorem
is equivalent to the Hamiltonian of the cases of Kovalevskaya or Goryachev. 

From  Corollary \ref{typ}
it follows immediately that in the family (\ref{gor}) of Goryachev 
we need to consider only the case $B_1=B_2=0$ which is in fact the case of 
Kovalevskaya, see the Introduction. 

Let us write the Hamiltonian (\ref{kov1}) of the case of Kovalevskaya in polar coordinates. 
We obtain
$$H=\gamma_1(r^2)(r^2d\varphi^2+dr^2)-\gamma_2(r^2)\cos\varphi$$
where $\gamma_2(r^2)\ne 0$ for all $0<r<+\infty$. Comparing  this with (\ref{H1})
 where 
$u\in (-\infty, +\infty)$, we see that no Hamiltonian  from our theorem
is equivalent to the Hamiltonian of the case of Kovalevskaya.

\hfill$\Box$

\begin{theo}
Assume  $b=1$. Then for any  $a>2$, 
$p\in (-\infty,-\frac{a}{2})\cup (-1,+\infty)$   and $-2<a<2$, 
$p\in(-\infty, -1)\cup(-\frac{a}{2}, +\infty)$ the Hamiltonian
$H_p$ (or $-H_p$) where $H_p$ has the form (\ref{H_p1}) and
$\varphi\in [0, 2\pi)$, $u\in (-\infty, +\infty)$ defines a  conservative system
on $S^2$ possessing an integral of fourth degree in momenta. This 
integral is nontrivial, i.e. there is no quadratic or linear integral.

These Hamiltonians for different values of parameters $a$ and $p$ are not equivalent. 
No Hamiltonian from this family is equivalent to the Hamiltonians of the cases of Kovalevskaya or Goryachev.
\label{th2}
\end{theo}

\noindent
{\em Proof.} The integrability follows immediately from Theorem
\ref{loc-th}. 

We will use  Proposition \ref{cen} to rewrite the corresponding Hamiltonians
in polar coordinates $\varphi, r=\exp y$:
$$H_{p}=\frac{\xi_a(r^2)+p+1}{\nu_a^2(r^2)}(r^2d\varphi^2 + dr^2)
 + \frac{\mu_a(r^2)}{\xi_a(r^2)+p+1}r\cos \varphi$$
$$=\frac{\xi_a(\tilde r^2)+p+1}{\nu_a^2(\tilde r^2)}
(\tilde r^2d\tilde \varphi^2 + d\tilde r^2) - 
\frac{\mu_a(\tilde r^2)}{\xi_a(\tilde r^2)+p+1}\tilde r\cos \tilde \varphi $$
where $\tilde r=\frac{1}{r}$, $\tilde \varphi=-\varphi$ and 
$$\xi_a(t)=\int_{t}^{1}{\mu_a(s)\nu_a^{-1}(s)ds}.$$

Now we must find the admissible values of the parameter $p$. 
We have $$p>\max_{z\ge 0} f(z)\quad \mbox{or}\quad p<\min_{z\ge 0} f(z)$$ 
 where $f(z)=z-\sqrt{1+az+z^2}$.
By computation we obtain $\min f(x)=f(0)=-1$, $\max f(x)=f(\infty)=-\frac{a}{2}$
 if $-2<a<2$ and $\min f(x)=f(\infty)=-\frac{a}{2}$, $\max f(x)=f(0)=-1$
if $a>2$. 
\medbreak
In order to prove all other statements of the theorem one must only repeat the arguments from the 
proof of Theorem \ref{th1}.

\hfill$\Box$

\section {The case of Kovalevskaya}

\begin{theo}
If $b=0$, $a=1$, $p=0$, then the Hamiltonian (\ref{H_p1}) where 
$\varphi\in [0, 2\pi)$ and $u\in [0, +\infty)$ defines a conservative
system on $S^2$, corresponding to the case of Kovalevskaya.
\label{th3}
\end{theo}

\noindent
{\em Proof.} In this case
 the Hamiltonian (\ref{H_p1})  has the form
\begin{equation}
H=\frac{1}{\sqrt{1+u^2}(\sqrt{1+u^2}+u)}\left(d\varphi^2+\frac{du^2}{u\sqrt{1+u^2}}\right)+
\frac{1}{2(\sqrt{1+u^2}+u)}\cos\varphi
\label{kov2}
\end{equation}
where $u\in [0, +\infty)$, $\varphi\in [0,2\pi)$.

Let us introduce  new variables
$$x=\Psi(u)\cos\varphi, \ y=\Psi(u)\sin\varphi, \ z=\pm\sqrt{1-\Psi(u)^2}$$
where $$\Psi(u)=\frac{1}{\sqrt{1+u^2}+u}=\sqrt{1+u^2}-u, \  0<\Psi(u)<1.$$

We get
$$\sqrt{1+u^2}=\frac{1}{2}\left(\Psi(u)+\frac{1}{\Psi(u)} \right), \ 
u=\frac{1}{2}\left(-\Psi(u)+\frac{1}{\Psi(u)} \right).$$

Let us compute
$$dx^2 + dy^2 + 2dz^2=\Psi'^2du^2+\Psi^2d\varphi^2+2\frac{\Psi^2\Psi'^2}{1-\Psi^2}du^2=\Psi^2
\left(d\varphi^2+\frac{du^2}{u\sqrt{1+u^2}}\right).$$
Then (\ref{kov2}) can be rewritten in the form
$$H=2\frac{\Psi^2}{\Psi^2 + 1}\left(d\varphi^2+\frac{du^2}{u\sqrt{1+u^2}}\right) +
\frac{1}{2}\Psi\cos\varphi=
2\frac{dx^2+dy^2+2dz^2}{2x^2+2y^2+z^2}+\frac{1}{2}x$$
where
$$x^2+y^2+z^2=1.$$

\hfill$\Box$

\medbreak
\noindent
{\bf Acknowledgement.} 
E. Selivanova would like to thank Professor Gerhard  Huisken and the Arbeitsbereich 
"Analysis" of the University of  T\"ubingen 
 for their hospitality.
\medbreak
\ \

\medbreak

\end{document}